\documentclass[a4pitaper,11pt]{amsart}
\usepackage{amsmath,amsthm,amssymb}
\usepackage[mathscr]{eucal}
 \usepackage{cite}
\usepackage{upgreek}
\usepackage[bookmarks,bookmarksnumbered]{hyperref}
\usepackage{enumerate}

\numberwithin{equation}{section}

\newcommand{\car}{\curvearrowright}

\theoremstyle{plain}
\newtheorem{main}{Theorem}

\newtheorem{theorem}{Theorem}[section]

\newtheorem{lemma}[theorem]{Lemma}
\newtheorem{proposition}[theorem]{Proposition}

\theoremstyle{definition}

\newtheorem{definition}[theorem]{Definition}

\newtheorem{notation}[theorem]{Notation}
\newtheorem{remark}[theorem]{Remark}

\begin{document}

\title[Orbit equivalence rigidity for product actions]
{Orbit equivalence rigidity for product actions}

\author[D. Drimbe]{Daniel Drimbe}
\address{Department of Mathematics, University of Regina, 3737 Wascana Pkwy, Regina, SK S4S 0A2, Canada.}
\email{daniel.drimbe@uregina.ca}
\thanks{The author was partially supported by PIMS fellowship.}

\begin{abstract} 
Let $\Gamma_1,\dots,\Gamma_n$ be hyperbolic, property (T) groups, for some $n\ge 1$.
We prove that if a product $\Gamma_1\times\dots\times\Gamma_n \car X_1\times\dots\times X_n$ of measure preserving actions is stably orbit equivalent to a measure preserving action $\Lambda\car Y$, then $\Lambda\car Y$ is induced from an action $\Lambda_0\car Y_0$ such that
there exists a direct product decomposition $\Lambda_0=\Lambda_1\times\dots\times\Lambda_n$ into $n$ infinite groups. Moreover, there exists a measure preserving action $\Lambda_i\car Y_i$ that is stably orbit equivalent to $\Gamma_i\car X_i$, for any $1\leq i\leq n$, and the product action $\Lambda_1\times\dots\times\Lambda_n\car Y_1\times\dots\times Y_n$ is isomorphic to $\Lambda_0\car Y_0$.

\end{abstract}

\maketitle

\section{Introduction}
An important topic in ergodic theory is the classification of probability measure preserving (pmp) actions up to orbit equivalence. Two pmp actions $\Gamma\car (X,\mu)$ and $\Lambda\car (Y,\nu)$ are called {\it orbit equivalent} (OE) if there exists a measure space isomorphism $\theta: (X,\mu)\to (Y,\nu)$ which preserves the orbits, i.e. $\theta(\Gamma x)=\Lambda \theta(x)$, for almost every $x\in X.$ 
The classification of actions up to OE is driven by the following fundamental question: what aspects of the group $\Gamma$ and of the action $\Gamma\car (X,\mu)$ are remembered by the {\it orbit equivalence relation} $\mathcal R_{\Gamma\car X}:=\{(x,y)\in X\times X| \Gamma x=\Gamma y\}$?

Equivalence relations $\mathcal R_{\Gamma\car X}$ tend to forget a lot of information about the groups and actions they are constructed from. This is best illustrated by H. Dye's theorem asserting that any two ergodic pmp actions of $\mathbb Z$ are OE \cite{Dy58}. D.S. Orstein and B. Weiss have extended this result to the class of countable amenable groups \cite{OW80} (see also \cite{CFW81} for a generalization). Consequently, pmp actions of amenable groups manifest a striking lack of rigidity: any algebraic property of the group (e.g. being finitely generated or torsion free) and any dynamical property of the action (e.g. being mixing or weakly mixing) is completely lost in the passage to equivalence relations. 

In the non-amenable case, the situation is radically different. More precisely, various properties of the group $\Gamma$ or of the action $\Gamma\car (X,\mu)$ can be recovered from the equivalence relation $\mathcal R_{\Gamma\car X}$. R. Zimmer's pioneering work led to such OE rigidity results for actions of higher rank lattices in semisimple Lie groups. In particular, he showed that if $m,n\ge 3$, then $SL_m(\mathbb Z)\car \mathbb T^m$ is OE to $SL_n(\mathbb Z)\car \mathbb T^n$ if and only if $m=n$ \cite{Zi84}. Remarkably, by building upon Zimmer's ideas, A. Furman showed that most pmp ergodic actions $\Gamma\car (X,\mu)$ of higher rank lattices, including $SL_n(\mathbb Z)\car \mathbb T^n$, for $n\ge 3$, are {\it OE superrigid}\cite{Fu99a,Fu99b}. Roughly speaking, this means that both the group $\Gamma$ and the action $\Gamma\car (X,\mu)$ are completely remembered by the equivalence relation $\mathcal R_{\Gamma \car X}$. By using his influential deformation/rigidity theory, S. Popa showed that any Bernoulli action of a property (T) group \cite{Po05} or of a product group \cite{Po06} is OE superrigid.
Subsequently, several impressive OE superrigidity results have been discovered in \cite{MS02,Ki06,Io08,PV08,Ki09,PS09,Io14,TD14,CK15,Dr15,GITD16}.

However, in general, one can only expect to recover certain aspects of an action $\Gamma\car (X,\mu)$ from its orbit equivalence relation $\mathcal R_{\Gamma\car X}$. For instance, D. Gaboriau showed that the rank of a free group $\mathbb F_{n}$ is an invariant of the OE relation $\mathcal R_{\mathbb F_n \car X}$ \cite{Ga99}. In \cite{MS02}, N. Monod and Y. Shalom obtained a series of OE rigidity results for actions of products $\Gamma=\Gamma_1\times\dots\times\Gamma_n$ of hyperbolic groups. In particular, they showed that if a product $\Gamma=\Gamma_1\times\dots\times\Gamma_m$ of non-elementary hyperbolic torsion-free groups is {\it measure equivalent} to a product $\Lambda=\Lambda_1\times\dots\times\Lambda_n$ of torsion-free groups (i.e. $\Gamma$ and $\Lambda$ admit stably orbit equivalent actions, see Definition \ref{D:ME}), then $m\ge n$, and if $m=n$ then, after a permutation of indices, $\Gamma_i$ is measure equivalent to $\Lambda_i$, for any $1\leq i\leq n$ (see \cite[Theorem 1.16]{MS02} and \cite[Theorem 3]{Sa09}). See the surveys \cite{Sh04,Po07,Fu09,Ga10,Va10,Io12,Io17} for an overview on orbit equivalence rigidity results and related topics.

The goal of this paper is to present a new type of rigidity phenomenon in orbit equivalence. Thus, we provide a large class of product actions whose orbit equivalence relation remember the product structure. Before stating the result, let us recall some terminology.

\begin{definition} Let $\Gamma\car (X,\mu)$ and $\Lambda\car (Y,\nu)$ be free ergodic pmp actions.
\begin{itemize}

\item The actions are {\it stably orbit equivalent} (SOE) {\it with index t} if there exists a measure space isomorphism $\theta: X_0\to Y_0$ between some non-null subsets $X_0\subset X$ and $Y_0\subset Y$ such that $\theta(X_0\cap \Gamma x)=Y_0\cap \Lambda \theta(x),$ for almost every $x\in X_0$. Moreover, $t$ is equal to $\nu(Y_0)/\mu(X_0).$

\item We say that $\Lambda\car Y$ is {\it induced from} an action $\Lambda_0\car Y_0$ of a subgroup $\Lambda_0<\Lambda$ if $Y_0\subset Y$ is a $\Lambda_0$-invariant measurable non-null subset such that $\nu(gY_0\cap Y_0)=0$, for all $g\in\Lambda\setminus\Lambda_0.$ In this case, note that $[\Lambda:\Lambda_0]=\nu(Y_0)^{-1}.$

\end{itemize}

\begin{remark}\label{induce}
We will use the following observation about induction. Let $\Lambda\car (Y,\nu)$ be a free ergodic pmp action.

\begin{itemize}

\item If $\Lambda\car Y$ is induced from $\Lambda_1\car Y_1$ and $\Lambda_1\car Y_1$ is induced from $\Lambda_0\car Y_0$, then $\Lambda\car Y$ is induced from $\Lambda_0\car Y_0$.

\item If $\Lambda\car Y$ is induced from $\Lambda_0\car Y_0$, then $\Lambda\car Y$ is SOE to $\Lambda_0\car Y_0$ with index $\nu(Y_0)^{-1}$.

\end{itemize}

\end{remark}

\end{definition}

\begin{main}\label{A2}\label{A}
Let $\Gamma_1,\dots,\Gamma_n$ be property (T), biexact, weakly amenable groups and denote $\Gamma=\Gamma_1\times ... \times\Gamma_n$. For every $1\leq i\leq n$, let $\Gamma_i\overset{}{\car} (X_i,\mu_i)$ be a free ergodic pmp action, and let $\Gamma\car (X,\mu)$ be the product action $\Gamma_1\times\dots\times\Gamma_n\car (X_1\times ... \times X_n, \mu_1\times ... \times \mu_n)$.\\ Let $\Lambda\car (Y,\nu)$ be a free ergodic pmp action of an icc group which is SOE to $\Gamma\car (X,\mu)$ with index $t$, for some $t>0$. 

Then $\Lambda\car Y$ is induced from an action $\Lambda_0\car Y_0$ and $\Lambda_0=\Lambda_1\times ... \times \Lambda_n$ decomposes into a direct product of $n$ infinite groups. Moreover, there exist a pmp action $\Lambda_i\car (Y_i,\nu_i)$ and a positive number $t_i$, for any $1\leq i\leq n$, with $t_1 \dots t_n=t/[\Lambda:\Lambda_0]$ such that:
\begin{enumerate}
\item $\Lambda_0\car Y_0$ is isomorphic to $\Lambda_1\times ... \times \Lambda_n\car Y_1\times ...\times Y_n$.

\item $\Lambda_i\car Y_i$ is SOE to $\Gamma_i\car X_i$ with index $t_i$, for all $1\leq i\leq n.$
\end{enumerate}

\end{main} 

Note that hyperbolic groups are weakly amenable and biexact, see for example the introductions of \cite{PV11,PV12} for an extensive discussion on these notions. 
Remark that the following classes of groups are property (T), biexact, weakly amenable groups:
\begin{itemize}
\item uniform lattices in $Sp(m,1)$ with $m\ge 2$ or any group in their measure equivalence class, 
\item Gromov's random groups with density satisfying $3^{-1}<d<2^{-1}$ \cite{Gr93, Zu03}. 
\end{itemize} 

A result due to Singer \cite{Si55} provides an approach to the study of orbit equivalence of actions using von Neumann algebras. In this context, we use a combination of techniques from S. Popa's deformation/rigidity theory to prove Theorem \ref{A2}.

Our main result can be seen as the orbit equivalence analog of a striking theorem of I. Chifan, R. de Santiago and T. Sinclair \cite[Theorem A]{CdSS15} in which they show that the group von Neumann algebra $L(\Gamma)$ of a product of non-elementary hyperbolic groups completely remembers the product structure of the generating group $\Gamma$.

Theorem \ref{A} complements the following remarkable result of N. Monod and Y. Shalom \cite[Theorem 1.9]{MS02}. First, we recall that a pmp action $\Lambda\car (Y,\nu)$ is {\it mildly mixing} if there are no nontrivial recurrent subsets: if a measurable subset $A\subset Y$ satisfies liminf$_{g\to\infty} \nu(gA\Delta A)=0$, then $\nu(A)\in\{0,1\}$. 

\begin{theorem}[\!\!\cite{MS02}]
Let $\Gamma=\Gamma_1\times\Gamma_2$ be a product of torsion-free non-elementary hyperbolic groups and let $\Gamma\overset{}{\car} X$ be a free irreducible pmp action. Assume $\Gamma\overset{}{\car} X$ is orbit equivalent to a mildly mixing action $\Lambda\overset{}{\car} Y$, where $\Lambda$ is a torsion-free group. 

Then there exists a group isomorphism $\delta:\Gamma\to\Lambda$ and a measure space isomorpshism $\theta:X\to Y$ such that $\theta(gx)=\delta(g)\theta(x)$, for all $g\in \Gamma$ and almost every $x\in X$.

\end{theorem}

Remark that Monod and Shalom's result applies to the class of {\it irreducible actions} $\Gamma_1\times\Gamma_2\overset{\sigma}{\car} X$, i.e. $\sigma_{|\Gamma_i}$ is ergodic for any $1\leq i\leq 2$, while Theorem \ref{A2} provides orbit equivalence rigidity results for product actions.

{\bf Comments on the proof of Theorem \ref{A2}.} 
We outline briefly and informally the proof of our main result. Assume for simplicity that $\Gamma=\Gamma_1\times\Gamma_2$ is a product of two hyperbolic, property (T) groups. Let $\Gamma_i\car X_i$ be a free ergodic pmp action, for any $1\leq i\leq 2.$ Denote by $\Gamma\car X$ the product action $\Gamma_1\times\Gamma_2\car X_1\times X_2$ and assume that it is OE to a free aperiodic action $\Lambda\car Y$, where $\Lambda$ is an icc group. By a result of Singer \cite{Si55} (see also \cite{FM75}) there exists an isomorphism of the groups measure space von Neumann algebras, $L^\infty(X)\rtimes\Gamma\cong L^\infty(Y)\rtimes\Lambda$, which identifies $L^\infty(X)$ and $L^\infty(Y).$ Hence, we assume $M:=L^\infty(X)\rtimes\Gamma=L^\infty(Y)\rtimes\Lambda$ and $L^\infty(X)=L^\infty(Y).$

Following \cite{PV09}, we define the comultiplication $*$-homomorphism $\Delta: M\to M\bar\otimes L(\Lambda)$ by letting $\Delta(av_g)=av_g\otimes v_g$, for all $a\in L^\infty(X), g\in\Lambda.$ We use the relative strong solidity property of hyperbolic groups (see Section \ref{S:rss}) obtained in the breakthrough work of S. Popa and S. Vaes \cite{PV11,PV12} and the property (T) assumption of $\Gamma_1$ to
 deduce that there exists $j\in\{1,2\}$ such that 
\begin{equation}\label{q1}
\Delta(L(\Gamma_1))\prec_{M\bar\otimes M} M \bar\otimes (L^\infty(X_j)\rtimes\Gamma_{j}).
\end{equation}
Here, $P\prec_N Q$ denotes that a corner of $P$ embeds into $Q$ inside $N$, in the sense of Popa \cite{Po03}.

A key step of the proof consists in applying A. Ioana's ultrapower technique from \cite{Io11} (see Theorem \ref{Th: ultrapower}) which allows us to deduce from \eqref{q1} the existence of a subgroup $\Sigma<\Lambda$ such that 
\begin{equation}\label{q2}
L^\infty(X)\rtimes\Gamma_1\prec_M L^\infty(X)\rtimes\Sigma \text{ and the centralizer of $\Sigma$ in $\Lambda$ is non-amenable.}
\end{equation}

By using \eqref{q2} together with the relative strong solidity of the groups $\Gamma_i$'s, we obtain that there exists a non-zero projection $e\in (L^\infty(X)\rtimes\Sigma)'\cap M$ such that 
\begin{equation}\label{q3}
L^\infty(X)\rtimes\Gamma_1\prec_M (L^\infty(X)\rtimes\Sigma)f \text{ and }
(L^\infty(X)\rtimes\Sigma)f \prec_M L^\infty(X)\rtimes\Gamma_1,
\end{equation}
for all $f\in (L^\infty(X)\rtimes\Sigma)'\cap M$ with $f\leq e.$

Next, by applying \cite[Proposition 3.1]{DHI16} to \eqref{q3} we derive that $\Sigma$ is measure equivalent to $\Gamma_1$, and hence, obtain that $\Sigma$ has property (T). Inspired by techniques from \cite{CdSS15} we further deduce that, up to replacing $\Sigma$ by a finite index subgroup, the subgroup generated by $\Sigma$ and its centralizer $C_\Lambda(\Sigma)$ has finite index in $\Lambda.$ Finally, by using once again \eqref{q3} and adapting some arguments from \cite{CdSS15}, we conclude that $\Lambda\car Y$ is induced from an action $\Lambda_0\car Y_0$ that admits a direct product decomposition $\Lambda_1\times\Lambda_2\car Y_1\times Y_2$ such that $\Lambda_i\car Y_i$ is stably orbit equivalent to $\Gamma_i\car X_i$, for any $1\leq i\leq 2.$ 

{\bf Acknowledgment.} I would like to thank Adrian Ioana for many comments that helped improve the exposition of the paper. I am also grateful to the referees for carefully reading the paper, for all their  suggestions and for pointing out an error in an earlier version of Theorem \ref{P: splits}.

\section{Preliminaries}

\subsection{Terminology} 
In this paper we consider {\it tracial von Neumann algebras} $(M,\tau)$, i.e. von Neumann algebras $M$ equipped with a faithful normal tracial state $\tau: M\to\mathbb C.$ This induces a norm on $M$ by the formula $\|x\|_2=\tau(x^*x)^{1/2},$ for all $x\in M$. We will always assume that $M$ is a {\it separable} von Neumann algebra, i.e. the $\|\cdot\|_2$-completion of $M$ denoted by $L^2(M)$ is separable as a Hilbert space.
We denote by $\mathcal U(M)$ the {\it unitary group} of $M$ and by $\mathcal Z(M)$ its {\it center}.

All inclusions $P\subset M$ of von Neumann algebras are assumed unital. We denote by $E_{P}:M\to P$ the unique $\tau$-preserving {\it conditional expectation} from $M$ onto $P$, by $P'\cap M=\{x\in M|xy=yx, \text{ for all } y\in P\}$ the {\it relative commutant} of $P$ in $M$ and by $\mathcal N_{M}(P)=\{u\in\mathcal U(M)|uPu^*=P\}$ the {\it normalizer} of $P$ in $M$.  We say that $P$ is {\it regular} in $M$ if the von Neumann algebra generated by $\mathcal N_M(P)$ equals $M$.
For two von Neumann subalgebras $P,Q\subset M$, we denote by $P\vee Q$ the von Neumann algebra generated by $P$ and $Q$. 

The {\it amplification} of a II$_1$ factor $(M,\tau)$ by a positive number $t$ is defined to be $M^t=p(\mathbb B(\ell^2(\mathbb Z))\bar\otimes M)p$, for a projection $p\in \mathbb B(\ell^2(\mathbb Z))\bar\otimes M$ satisfying $($Tr$\otimes\tau)(p)=t$. Here Tr denotes the usual trace on $\mathbb B(\ell^2(\mathbb Z))$. Since $M$ is a II$_1$ factor, $M^t$ is well defined. Note that if $M=P_1\bar\otimes P_2$, for some II$_1$ factors $P_1$ and $P_2$, then there exists a natural identification $M=P_1^t\bar\otimes P_2^{1/t}$, for every $t>0.$

Let $\Gamma\overset{\sigma}{\car} A$ be a trace preserving action of a countable group $\Gamma$ on a tracial von Neumann algebra $(A,\tau)$. For a subgroup $\Sigma<\Gamma$, we denote by $A^\Sigma=\{a\in A|\sigma_g(a)=a, \text{ for all } g\in\Sigma\}$, the subalgebra of elements of $A$ fixed by $\Sigma.$

For a countable group $\Gamma$ and for two subsets $S,T\subset \Gamma$, we denote by $\langle S \rangle$ the group generated by $S$, and by $C_S(T)=\{g\in S|gh=hg, \text{ for all } h\in T\}$ the {\it centralizer of $T$ in $S$}. 

For an abelian von Neumann algebra $A=L^\infty(X)$ and a measurable subset $X_0\subset X$, we denote by $1_{X_0}$ the associated projection in $A$.

\subsection {Intertwining-by-bimodules} We next recall from  \cite [Theorem 2.1 and Corollary 2.3]{Po03} the powerful {\it intertwining-by-bimodules} technique of S. Popa.

\begin {theorem}[\!\!\cite{Po03}]\label{corner} Let $(M,\tau)$ be a tracial von Neumann algebra and $P\subset pMp, Q\subset qMq$ be von Neumann subalgebras. Let $\mathcal U\subset\mathcal U(P)$ be a subgroup such that $\mathcal U''=P$.

Then the following are equivalent:

\begin{enumerate}

\item There exist projections $p_0\in P, q_0\in Q$, a $*$-homomorphism $\theta:p_0Pp_0\rightarrow q_0Qq_0$  and a non-zero partial isometry $v\in q_0Mp_0$ such that $\theta(x)v=vx$, for all $x\in p_0Pp_0$.

\item There is no sequence $(u_n)_n\subset\mathcal U$ satisfying $\|E_Q(xu_ny)\|_2\rightarrow 0$, for all $x,y\in M$.

\end{enumerate}
\end{theorem}

\begin{notation} Throughout the paper we will use the following notation.
\begin{itemize}
\item If one of the equivalent conditions of Theorem \ref{corner} holds true, we write $P\prec_{M}Q$, and say that {\it a corner of $P$ embeds into $Q$ inside $M$.}
\item If $Pp'\prec_{M}Q$ for any non-zero projection $p'\in P'\cap pMp$ (equivalently, for any non-zero projection $p'\in \mathcal N_{pMp}(P)'\cap pMp$, by Lemma \ref{Lemma 2.4, DHI16}(1)), then we write $P\prec^{s}_{M}Q$.
\item If $P\prec_{M}Qq'$ for any non-zero projection $q'\in Q'\cap qMq$ (equivalently, for any non-zero projection $q'\in \mathcal N_{qMq}(Q)'\cap qMq$, by Lemma \ref{Lemma 2.4, DHI16}(2)), then we write $P\prec^{s'}_{M}Q$.

\end{itemize}
\end{notation}

We continue with some results containing several elementary facts regarding Popa's intertwining-by-bimodules technique.

\begin{lemma}\label{extension} 
Let $(M,\tau)$ be a tracial von Neumann algebra and let $N\subset M$ be a von Neumann subalgebra. Let $P\subset pNp$ and $Q\subset qNq$ be von Neumann subalgebras such that $P\prec_N^{s'} Q$.

Then $P\prec_M^{s'} Q.$
\end{lemma}

Note that Lemma \ref{extension} also holds true if we replace the symbol $\prec^{s'}$ by $\prec^s$ as shown in \cite[Remark 2.2]{DHI16}.

{\it Proof.} Let $q'\in Q'\cap qMq$ be a non-zero projection. Let $q''\in Q'\cap qNq$ be the support projection of $E_N(q').$ The assumption implies that $P\prec_{N}Qq''$, hence there exist projections $p_0\in P, q_0\in Q$, a $*$-homomorphism $\theta: p_0Pp_0\to q_0Qq_0q''$ and a non-zero partial isometry $v\in q_0q''Np_0$ such that $\theta(x)v=vx$, for all $x\in p_0Pp_0.$ Since $q'\leq q''$, we let $\tilde\theta: p_0Pp_0\to q_0Qq_0q'$ be the $*$-homomorphism defined by $\tilde\theta(x)=\theta(x)q',$ for all $x\in p_0Pp_0$. By denoting $\tilde v=q'v$, we have $\tilde\theta (x)\tilde v=\tilde v x$, for all $x\in p_0Pp_0$. Note that $\tilde v$ is non-zero since $E_{N}(q')v\neq 0$. Finally, by replacing $\tilde v$ by the partial isometry from its polar decomposition, we obtain that $P\prec_M Qq'.$
\hfill$\blacksquare$

\begin{lemma}\label{transitivity}
Let $(M,\tau)$ be a tracial von Neumann algebra and let $P\subset pMp$, $Q\subset qMq$, $R\subset rMr$ be von Neumann subalgebras. Then the following hold:

\begin{enumerate}

\item  Assume that $P\prec_M Q$ and $Q\prec_M^s R$. Then $P\prec_M R$, \emph{\cite[Lemma 3.7]{Va08}}.

\item Assume $P\prec_M^{s'} Q$ and $Q\prec_M R$. Then $P\prec_M R.$

\end{enumerate}

\end{lemma}

{\it Proof.} We will prove only the second statement.

(2) The assumption $Q\prec_M R$ implies that there exist projections $q_0\in Q, r_0\in R$, a non-zero partial isometry $v\in r_0Mq_0$ and a $*$-homomorphism $\psi: q_0Qq_0\to r_0Rr_0$ such that $\psi(x)v=vx$, for all $x\in q_0Qq_0.$ Let $q'$ be the support projection of $v^*v$ and note that $q'\in Q'\cap qMq$ and $q'\leq q_0$. Without loss of generality, we can assume that the support projection of $E_{R}(vv^*)$ equals $r_0.$
The assumption implies that $P\prec_M Qq'$, hence there exist projections $p_0\in P, q_1\in Q$, a non-zero partial isometry $w\in q_1q'Mp_0$ and a $*$-homomorphism $\phi: p_0Pp_0\to q_1Qq_1q'$ such that 
$\phi(x)w=wx$, for all $x\in p_0Pp_0$. We can assume that the support projection of $E_{Qq'}(ww^*)$ equals $q'q_1.$

Since $q_0Qq_0q'=Qq'$, we can define the $*$-homomorphism $\psi':Qq'\to r_0Rr_0$ by letting $\psi'(xq')=\psi(x)$, for all $x\in q_0Qq_0.$ Notice that $\psi'$ is well defined because if $xq'=0$, for some $x\in q_0Qq_0$, then $\psi (x)vq'=0.$ This shows that $\psi (x)E_{R}(vv^*)=0$, which implies $\psi(x)=0$. Moreover, it is clear that $\psi'(xq')vq'=vq'(xq')$, for all $x\in q_0Qq_0.$ We continue by noticing that $vq'w\neq 0$. Indeed, by assuming the contrary, we obtain that $q'E_{Qq'}(ww^*)=0,$ which implies that $q'q_1=0$, contradiction.

Finally, remark that the $*$-homomorphism $\psi' \circ \phi : p_0Pp_0\to r_0Rr_0$ satisfies $(\psi'\circ\phi)(x)vq'w=vq'wx$, for all $x\in p_0Pp_0$. By replacing $vq'w$ by the partial isometry from its polar decomposition, we deduce that $P\prec_M R.$ 
\hfill$\blacksquare$

We will use repeatedly the following results from \cite{DHI16} and \cite{Va08} and we record them in the following combined lemma for reader's convenience.

\begin{lemma}\label{Lemma 2.4, DHI16}\label{Lemma 3.5, Va08}\label{dhi}\label{va}
Let $(M,\tau)$ be a tracial von Neumann algebra and let $P\subset pMp$, $Q\subset qMq$ be von Neumann subalgebras. Then the following hold:

\begin{enumerate}
\item If $Pz\prec_M Q$, for any non-zero projection $z\in\mathcal N_{pMp}(P)'\cap pMp$, then $P\prec^s_M Q$, \emph{\cite[Lemma 2.4(2)]{DHI16}}.

\item If $P\prec_M Q$, then there exists a non-zero projection $z\in\mathcal N_{qMq}(Q)'\cap qMq$ such that $P\prec^{s'}_M Qz$, \emph{\cite[Lemma 2.4(4)]{DHI16}}.

\item If $P\prec_M Q$, then $Q'\cap qMq\prec_M P'\cap pMp$, \emph{\cite[Lemma 3.5]{Va08}}.

\end{enumerate}

\end{lemma}

If we consider two commuting subalgebras $P_1$ and $P_2$ of $M$ that both embed into a subalgebra $Q$ of $M$ in the sense of Popa (see Theorem \ref{corner}), we would like to obtain that $P_1\vee P_2$ embeds into $Q$. In general this is not true. For instance, consider the crossed product $M=(P\bar\otimes P)\rtimes_\sigma \mathbb Z/2\mathbb Z$ with the period two automorphism $\sigma(a\otimes b)=b\otimes a$, where $P$ is a tracial diffuse von Neumann algebra. Then, we have that $P\otimes 1\prec_M  P\otimes 1$ and $1\otimes P\prec_M P\otimes 1$, but $P\bar\otimes P\nprec_M P\otimes 1.$
The following lemma essentially shows that in the presence of a weak regularity condition on $Q$, the result is true. Although its proof is easy, this lemma plays an important role in the paper.

\begin{lemma}\label{L: full}
Let $M$ be a von Neumann algebra and let $P_1, P_2\subset M$ and $Q\subset qMq$ be von Neumans subalgebras such that $\mathcal {N}_{qMq}(Q)' \cap qMq=\mathbb C q$. Suppose there exist commuting subalgebras $\tilde P_0,\tilde P_1, \tilde P_2\subset M$ such that $P_1\subset \tilde P_1$, $P_2\subset \tilde P_2$ and $\tilde P_0\vee \tilde P_1\vee \tilde P_2=M$. 

If $P_i\prec_M Q$, for any $i\in \{1,2\},$  then $P_1\vee P_2\prec_M Q.$
\end{lemma}

{\it Proof.}
The assumption implies that there exist projections $p\in P_1, q_0\in Q$, a $*$-homomorphism $\varphi: pP_1p\to q_0Qq_0$ and a non-zero partial isometry $v\in q_0Mp$ such that $\varphi(x)v=vx$, for every $x\in pP_1p.$ 

We aim to show that $(pP_1p)\vee (P_2p)\prec_M Q.$ By supposing the contrary, there exist two sequences of unitaries $(u_n)_n\subset\mathcal U(pP_1p)$ and $(v_n)_n\subset\mathcal U(P_2)$, such that 
$$
\|E_Q(xu_n(v_np)y)\|_2\to 0, \text{ for all }x,y\in M.
$$
Since $M=\tilde P_0\vee \tilde P_1\vee \tilde P_2,$ we obtain that
\begin{equation}\label{eq: info1}
\|E_Q(xu_nyv_nz)\|_2\to 0, \text{ for all }x,y,z\in M.
\end{equation}

Using that $\varphi(u_n)v=vu_n$, for any $n\ge 1$, we have $\|E_{Q}(v u_n y v_nz\|_2 = \|E_{Q}(vy v_n z)\|_2$, for all $y,z\in M$. Combining this last remark with relation \eqref{eq: info1}, we obtain that $\|E_{Q}(ay v_n z)\|_2\to 0$, for all $y,z\in M$, where $a:=vv^*\in qMq.$ Since $wE_Q(x)w^*=E_Q(wxw^*)$, for any $x\in M$ and $ w\in\mathcal N_{qMq}(Q)$, it follows that 
\begin{equation}\label{N}
\|E_{Q}((waw^*)y v_n z)\|_2 \to 0, \text{ for all } y,z \in M \text{ and  }w\in\mathcal N_{qMq}(Q)
\end{equation}

We continue by noticing that for any two projections $p_1$ and $p_2$ in $M$, we have $p_1\vee p_2=s(p_1+p_2)$. Here we denote by $s(b)$ the support projection of an an element $b\in M.$ 
Moreover, by using Borel functional calculus, there exist a sequence $(c_n)_n\subset M$ such that $(p_1+p_2)c_n$ converges to $s(p_1+p_2)$ in the $\|\cdot\|_2$-norm. Thus, it follows by \eqref{N} that $\|E_{Q}((w_1aw_1^*+w_2aw_2^*) y v_n z)\|_2\to 0$, and hence, $\|E_{Q}((w_1aw_1^*\vee w_2aw_2^*) y v_n z)\|_2\to 0$, for any $w_1,w_2 \in\mathcal N_{qMq}(Q)$ and $y,z\in M$.

Therefore, we obtain by induction that
$$
\|E_{Q}((w_1aw_1^*\vee...\vee w_naw_n^*) yv_n z)\|_2\to 0,
$$
for any $w_1,...,w_n\in\mathcal N_{qMq}(Q)$ and $y,z\in M$. \\
Finally, remark that $\vee_{w\in\mathcal N_{qMq}(Q)}waw^*=q$ since $\mathcal {N}_{qMq}(Q)' \cap qMq=\mathbb C q$. Hence, it follows that $ \|E_Q(y v_n z)\|_2\to 0$, for all $y,z\in M$. This shows that $P_2\nprec_M Q$, contradiction.
\hfill$\blacksquare$

The next lemma is well known. We include only a sketch of the proof since the result follows, for example, from \cite[Lemma 1]{HPV11} (see also \cite[Theorem 6.2]{Po01}).

\begin{lemma}\label{L: rigid}
Let $M$ and $M_0$ be some tracial von Neumann algebras and let $A\subset M$ be an abelian subalgebra. Denote $\mathcal M=M_0\bar\otimes M.$ 

If $P\subset p\mathcal M p$ is a property (T) subalgebra such that $P\prec_{\mathcal M} M_0\bar\otimes A$, then $P\prec_{\mathcal M} M_0.$ 
\end{lemma}

{\it Proof.}
Since $A$ is abelian, take an increasing sequence of finite dimensional abelian algebras $A_n$, $n\ge 1.$ Using the fact that $P$ has property (T) and $P\prec_{\mathcal M} (\cup_n M_0\bar\otimes A_n)''$, we get that there exists an integer $n_0$ such that $P\prec_{\mathcal M} M_0\bar\otimes A_{n_0}.$ This proves the claim since $A_{n_0}$ is finite dimensional. 
\hfill$\blacksquare$

Property (T) for von Neumann algebras was defined by Connes and Jones in \cite{CJ85}. Note that a countable group $\Gamma$ has property (T) if and only if $L(\Gamma)$ has property (T) \cite{CJ85, Po01}.

The following result is a tool that was discovered in \cite[Proposition 3.1]{DHI16} which 
allows us to derive from some intertwining relations that certain groups are {measure equivalent} in the sense of M. Gromov. We first recall the notion of measure equivalence \cite{Gr91,Fu99a}.
\begin{definition}\label{D:ME}
Two countable groups $\Gamma$ and $\Lambda$ are called {\it measure equivalent} if there exist free ergodic pmp actions $\Gamma\car (X,\mu)$ and $\Lambda\car (Y,\nu)$ which are stably orbit equivalent.
\end{definition}

\begin{proposition}[\!\!\cite{DHI16}]\label{ME}
Let $M$ be a II$_1$ factor and let $\Gamma\car (X,\mu)$ and $\Lambda\car (Y,\nu)$ be free ergodic pmp actions such that $pMp=L^\infty(X)\rtimes\Gamma$ and $qMq=L^\infty(Y)\rtimes\Lambda$, for some projections $p,q\in M$. Suppose that $\Gamma=\Gamma_1\times\Gamma_2$ and $L^\infty(X)\prec_M L^\infty(Y)$. Assume that there exists a subgroup $\Sigma<\Lambda$ such that the following hold:
\begin{itemize}
\item $L^\infty(X)\rtimes\Gamma_1\prec_M L^\infty(Y)\rtimes\Sigma$, and
\item $L^\infty(Y)\rtimes\Sigma \prec_M^s L^\infty(X)\rtimes\Gamma_1.$
\end{itemize}
Then $\Sigma$ is measure equivalent to $\Gamma_1$.
\end{proposition}

\subsection{Finite index inclusions of von Neumann algebras}

The {\it Jones index} for an inclusion $P\subset M$ of II$_1$ factors is the dimension of $L^2(M)$ as a left $P$-module \cite{Jo81}. M. Pimsner and S. Popa defined a probabilistic notion of index for an inclusion $P\subset M$ of arbitrary von Neumann algebras with conditional expectation, which in the case of inclusions of II$_1$ factors coincides with Jones' index \cite[Theorem 2.2]{PP86}.
Namely, the inclusion $P\subset M$ of tracial von Neumann algebras is said to have {\it probabilistic index} $[M:P]=\lambda^{-1}$, where
$$
\lambda= \text{inf}\{\|E_P(x)\|^2_2{\|x\|_2^{-2}}|x\in M_+, x\neq 0\}.
$$
Here we use the convention that $\frac{1}{0}=\infty.$

We continue with recording several basic facts concerning finite index inclusions of von Neumann algebras. For a proof of the next lemma, see \cite[Lemma 2.4]{CIK13}, for example.

\begin{lemma}[\!\!{\cite[Lemma 2.3]{PP86}}]\label{fi1}
Let $N\subset M$ be tracial von Neumann algebras satisfying $[M:N]<\infty$. Then the following hold:

\begin{enumerate}
\item for every projection $p\in N$, we have
 $[pMp:pNp]<\infty$.
\item $M\prec_M^s N$.

\end{enumerate}

\end{lemma}

\begin{lemma}\label{fi2}
Let $N\subset M$ be tracial von Neumann algebras satisfying $[M:N]<\infty$
and assume that $\mathcal Z(N)$ is completely atomic.

Then $qMq\prec_{qMq} Nq$, for every projection $q\in N'\cap M.$
Moreover, $qMq\prec^{s'}_{qMq} Nq$, for every projection $q\in N'\cap M$.
\end{lemma}

{\it Proof.} Let $q\in N'\cap M$ be a non-zero projection. Let $z\in\mathcal Z(N)$ be a projection such that $Nz$ is a factor and $qz\neq 0$. Lemma \ref{fi1}(1) implies that $[zMz:Nz]<\infty$. Since $Nz$ is a factor, by using \cite[Proposition 2.3(3)]{CdSS17}, we get that $[qzMqz:Nqz]<\infty$, hence we obtain $qzMqz\prec_{qzMqz} Nqz$ by Lemma \ref{fi1}(2). Therefore, $qMq\prec_{qMq} Nq$. Note that the moreover part follows from the first part. 
\hfill$\blacksquare$

The following lemma goes back to \cite{Jo81} and extends the results from \cite[Examples 2.3.2 and 2.3.3]{Jo81}.

\begin{lemma}\label{fi3}
Let $\Gamma\overset{\sigma}{\car} (A,\tau)$ be a trace preserving action and denote $M=A\rtimes\Gamma$. Let $(\Delta_1,\Delta_2)$ and $(\Omega_1,\Omega_2)$ be two pairs of commuting subgroups of $\Gamma$ such that $\Delta_1\subset \Omega_1$ and $\Omega_2\subset \Delta_2$ are finite index inclusions.

Then the inclusion $A^{\Omega_1}\rtimes\Omega_2\subset A^{\Delta_1}\rtimes\Delta_2 $ has finite index.
\end{lemma}

{\it Proof.} Note that it is enough to show the following two statements:
\begin{enumerate}
\item $A^{\Omega_1}\rtimes\Omega_2 \subset A^{\Delta_1}\rtimes\Omega_2$ has finite index. 

\item $A^{\Delta_1}\rtimes\Omega_2 \subset A^{\Delta_1}\rtimes\Delta_2$ has finite index.
\end{enumerate}

(1) Denote by $\{u_g\}_{g\in\Gamma}$ the canonical unitaries which implement the action $\Gamma\car A$. Let $g_1,\dots,g_n\in \Omega_1$ such that we have the partition $\Omega_1=\underset{1\leq i\leq n}{\sqcup}g_i\Delta_1.$ One can check that the map $E_0:A^{\Delta_1}\to A^{\Omega_1}$ defined by $E_0(a)=\frac{1}{n}\sum_{i=1}^n \sigma_{g_i}(a)$ is the unique conditional expectation from $A^{\Delta_1}$ to $A^{\Omega_1}$.  Denote by $E:A^{\Delta_1}\rtimes\Omega_2\to A^{\Omega_1}\rtimes\Omega_2$ the conditional expectation from $A^{\Delta_1}\rtimes\Omega_2$ to $A^{\Omega_1}\rtimes\Omega_2$. Since $\Omega_1$ and $\Omega_2$ commute, note that for any $x=\sum_{g\in\Omega_2}x_gu_g\in A^{\Delta_1}\rtimes\Omega_2$, we have $E(x)=\sum_{g\in\Omega_2}E_0(x_g)u_g=\frac{1}{n}\sum_{i=1}^n u_{g_i} x  u_{g_i}^* $. Hence, for any $x\in (A^{\Delta_1}\rtimes\Omega_2)_+$, we have $\|E(x)\|_2^2\ge \frac{1}{n^2}\sum_{i=1}^n \tau (u_{g_i} x^*  u_{g_i}^* u_{g_i} x  u_{g_i}^*)=\frac{1}{n}\|x\|_2^2.$ This ends the first part.

(2) Denote $P=A^{\Delta_1}\rtimes\Omega_2$ and $N= A^{\Delta_1}\rtimes\Delta_2$. Denote by $e_P:L^2(N)\to L^2(P)$ the orthogonal projection onto $L^2(P)$ and note that the basic construction $\langle N,e_P\rangle $ is isomorphic to $(A^{\Delta_1}\bar\otimes \ell^\infty(\Delta_2/\Omega_2))\rtimes \Delta_2$ (see, for example,\cite[Lemma 2.5]{Be14}). Hence, $\langle N,e_P\rangle $ is tracial, since $[\Delta_2:\Omega_2]<\infty$. Therefore, there exists a normal conditional expectation $E:\langle N,e_P\rangle \to N$, which implies that $P\subset N$ is a finite index inclusion.
\hfill$\blacksquare$

\begin{lemma}\label{atomic}
Let $\Gamma\overset{\sigma}{\car} A$ be a trace preserving ergodic action on an abelian von Neumann algebra $(A,\tau)$ and let $\Sigma<\Gamma$ be a finite index subgroup.

Then, $A^\Sigma$ is completely atomic.

\end{lemma}

{\it Proof.} Let $n=[\Gamma:\Sigma]$ and take $g_1,\dots,g_n\in\Gamma$ such that we have the partition $\Gamma=\sqcup_{1\leq i \leq n}g_i\Sigma$. Assume by contrary that $A^\Sigma$ is not completely atomic. Then there exists a non-zero projection $p\in A^\Sigma$ such that $\tau(p)\leq 1/(2n)$. Note that $a:=\sum_{i=1}^n \sigma_{g_i}(p)$ belongs to $A^\Gamma=\mathbb C$ and $\tau(a)\leq 1/2$. This leads to a contradiction. 
\hfill$\blacksquare$

\subsection{Relative amenability of subalgebras}
A tracial von Neumann algebra $(M,\tau)$ is {\it amenable} if there exists a positive linear functional $\Phi:\mathbb B(L^2(M))\to\mathbb C$ such that $\Phi_{|M}=\tau$ and $\Phi$ is $M$-{\it central}, meaning $\Phi(xT)=\Phi(Tx),$ for all $x\in M$ and $T\in \mathbb B(L^2(M))$.

A very useful relative version of this notion has been introduced by N.
Ozawa and S. Popa in \cite{OP07}. Let $(M,\tau)$ be a tracial von Neumann algebra. Let $p\in M$ be a projection and $P\subset pMp,Q\subset M$ be von Neumann subalgebras. Following \cite[Definition 2.2]{OP07}, we say that $P\subset pMp$ is {\it amenable relative to $Q$ inside $M$} if there exists a positive linear functional $\Phi:p\langle M,e_Q\rangle p\to\mathbb C$ such that $\Phi_{|pMp}=\tau$ and $\Phi$ is $P$-central. 
Note that $P$ is amenable relative to $\mathbb C$ inside $M$ if and only if $P$ is amenable. 

We next record the following useful results:

\begin{lemma}[\!\!\cite{Be14,BV12}]\label{union}
Let $\Gamma\overset{\sigma}{\car} (X,\mu)$ be a pmp action and denote $A=L^\infty(X)$ and $ M=L^\infty(X)\rtimes\Gamma$.
Let $p\in A$ be non-zero projection and $\Sigma<\Gamma$ a subgroup. Let $\mathcal G$ be a subgroup of $\mathcal N_{pMp}(Ap)$.

\begin{enumerate}

\item If $\mathcal G''\prec_M A\rtimes\Sigma$, then $(Ap\cup \mathcal G)''\prec_M A\rtimes\Sigma$,
\emph{\cite[Lemma 2.3]{BV12}}.

\item Assume that $\sigma$ is free and $(\mathcal G q)''$ is amenable relative to $A\rtimes\Sigma$, for some non-zero projection $q\in \mathcal G'\cap pMp.$ Then, $(Ap\cup\mathcal G)''q_0$ is amenable relative to $A\rtimes\Sigma$, for some non-zero projection $q_0\in (Ap\cup\mathcal G)'\cap pMp$, \emph{\cite[Lemma 2.11]{Be14}}.

\end{enumerate}

\end{lemma}

\begin{proposition}[\!\!\cite{PV11,DHI16}]\label{intersection}
Let $(M,\tau)$ be a tracial von Neumann algebra and $Q_1, Q_2\subset M$ von Neumann subalgebras which form a commuting square, i.e. $E_{Q_1}\circ E_{Q_2}=E_{Q_2}\circ E_{Q_1}$. Assume that $Q_1$ is regular in $M$.
Let $P\subset pMp$ be a von Neumann subalgebra. Then the following hold:

\begin{enumerate} 

\item If $P\prec_M^s Q_1$ and $P\prec_M^s Q_2$, then $P\prec_{M}^sQ_1\cap Q_2$, \emph{\cite[Lemma 2.8(2)]{DHI16}}.

\item If $P$ is amenable relative to $Q_1$ and $Q_2$, then $P$ is amenable relative to $Q_1\cap Q_2$, \emph{\cite[Proposition 2.7]{PV11}}.

\end{enumerate} 
\end{proposition}

\subsection{Relatively strongly solid groups}\label{S:rss} 
Following \cite[Definition 2.7]{CIK13}, a countable group $\Gamma$ is said to be {\it relatively strongly solid} and write $\Gamma\in \mathcal C _{rss}$ if for any trace preserving action $\Gamma\car B$ the following holds: if $M=B\rtimes\Gamma$ and $A\subset pMp$ is a von Neumann subalgebra which is amenable relative to $B$ inside $M$, then either $A\prec_M B$ or the normalizer $\mathcal N_{pMp}(A)''$ is amenable relative to $B$ inside $M$.\\
In their breakthrough work \cite{PV11,PV12}, S. Popa and S. Vaes proved that non-elementary hyperbolic groups belong to $\mathcal C_{rss}.$ More generally, [PV12, Theorem 1.4] shows that all weakly amenable, biexact groups are relatively strongly solid.

The following consequence of belonging to $\mathcal C_{\text{rss}}$ will be useful (see \cite[Lemma 5.2]{KV15}).

\begin{lemma}[\!\!{\cite{KV15}}]\label{L: rss}
Let $\Gamma\car Q$ be a trace preserving action of a group $\Gamma$ that belongs to the class $\mathcal C_{\text{rss}}$, and let $M=Q\rtimes\Gamma$. 
Let $P_1, P_2\subset pMp$ be commuting von Neumann subalgebras.

Then either $P_1\prec_{M}Q$ or $P_2$ is amenable relative to $Q$.
\end{lemma}

\section{From tensor decompositions to decompositions of actions}

For proving Theorem \ref{A2} we need the following result, which provides sufficient conditions at the von Neumann algebra level for a pmp action to admit a non-trivial direct product decomposition (see also \cite[Theorem 3.1]{Dr19}).
The factor setting is essential for the result and its proof is based on arguments from \cite[Theorem 4.14]{CdSS15} (see also \cite[Theorem 6.1]{DHI16} and \cite[Theorem 4.7]{CdSS17}).

\begin{theorem}\label{P: splits}
Let $\Lambda\car (Y,\nu)$ be a free ergodic pmp action of an icc group $\Lambda$. Let $M=L^\infty(Y)\rtimes\Lambda$ and assume that $M=P_1\bar\otimes P_2$ for some II$_1$ factors $P_1$ and $P_2$.

Suppose that there exist infinite commuting subgroups $\Sigma_1,\Sigma_2<\Lambda$ with $[\Lambda:\Sigma_1\Sigma_2]<\infty$ such that 
$$
L(\Sigma_1)e\prec_M P_1   \text{ and }    L(\Sigma_2)e\prec^s_M P_2,
$$
for some non-zero projection $e\in L^\infty(Y)^{\Sigma_1\Sigma_2}$.

Then there exist infinite commuting subgroups $\Lambda_1,\Lambda_2<\Lambda$ such that $\Lambda\car Y$ is induced from an action $\Lambda_1\Lambda_2\car Y_0.$ 

Moreover, there exist
a decomposition $1_{Y_0}M1_{Y_0}=P_1^{t_1}\bar\otimes P_2^{t_2}$, for some $t_1,t_2>0$ with $\nu(Y_0)=t_1t_2$ and a unitary $u\in\mathcal U(1_{Y_0}M1_{Y_0})$ such that
$$
P_1^{t_1}=u (L^\infty(Y)^{\Lambda_2}\rtimes\Lambda_1)1_{Y_0}u^* \text{ and } 
P_2^{t_2}=u (L^\infty(Y)^{\Lambda_1}\rtimes\Lambda_2)1_{Y_0}u^*.
$$
In particular, there exists a pmp action $\Lambda_i\car (Y_i,\nu_i)$, for any $i\in\{1,2\}$, such that $\Lambda_1\Lambda_2 \car Y_0$ is isomorphic to the product action $\Lambda_1\times\Lambda_2 \car Y_1\times Y_2.$
\end{theorem}

Before proceeding with the proof, we introduce some terminology and recall a well known lemma. Let $\Sigma<\Lambda$ be a subgroup. Following \cite{CdSS15}, we denote by $\mathcal O_\Sigma (g)=\{hgh^{-1}|h\in\Sigma\}$ the orbit of $g\in \Lambda$ under the conjugation action of $\Sigma.$ Note that $\mathcal O_\Sigma (g_1g_2)\subset \mathcal O_\Sigma (g_1) \mathcal O_\Sigma (g_2),$ thus $|\mathcal O_\Sigma (g_1g_2)|\leq |\mathcal O_\Sigma (g_1)||\mathcal O_\Sigma (g_2)|$. This implies that the set $\Delta=\{g\in\Lambda|\mathcal O_\Sigma(g) \text{ is finite}\}$ is a subgroup of $\Lambda$. Note also that $\Delta$ is normalized by $\Sigma$. Moreover, one can check that if $\Lambda\car (Y,\nu)$ is a pmp action, then $L(\Sigma)'\cap (L^\infty(Y)\rtimes\Lambda)\subset L^\infty(Y)\rtimes \Delta.$

\begin{lemma}\label{not diffuse}
Let $(M,\tau) $ be a tracial von Neumann algebra and assume $M=P_1\bar\otimes P_2$ for some von Neumann algebras $P_1$ and $P_2$. Let $A\subset P_1$ be a von Neumann subalgebra such that $A\prec_M P_2$.

Then, $A$ is not diffuse.
\end{lemma}

{\it Proof.} Assume by contradiction that $A$ is diffuse. Therefore, there exists a sequence $(u_n)_n\subset A$ of unitaries such that $\tau(u_nx)\to 0$ for all $x\in A$. This implies that 
\begin{equation}\label{nd}
\|E_{P_2}(xu_ny)\|_2\to 0, \text{ for all } x,y\in M.
\end{equation}
Indeed, note that since $M=P_1\bar\otimes P_2$, we can assume $x=1$ and $y=y_1\otimes y_2$ with $y_1\in P_1$ and $y_2\in P_2$. Hence, $\|E_{P_2}(xu_ny)\|_2=\tau(u_ny_1)\|y_2\|_2\to 0$. Therefore, \eqref{nd} is true, which implies the contradiction $A\nprec_M P_2.$
\hfill$\blacksquare$

{\it Proof of Theorem \ref{P: splits}.} Let $A=L^\infty(Y)$ and denote $\Delta_2=\{g\in\Lambda| \mathcal O{_{\Sigma_1}}(g) \text{ is finite}\}.$ Since $\Sigma_2\subset \Delta_2$ and $[\Lambda:\Sigma_1\Sigma_2]<\infty$, we obtain that there exist $g_1,...,g_n\in \Delta_2$ such that $\Delta_2 \Sigma_1=\cup_{i=1}^n g_i \Sigma_2\Sigma_1.$ Since $[\Sigma_1:C_{\Sigma_1}(g_i)]<\infty$ for every $1\leq i\leq n$, we obtain that $\Delta_1:=\cap_{i=1}^n C_{\Sigma_1}(g_i)$ is a finite index subgroup of $\Sigma_1.$ Note that $\Sigma_1$ is icc, since $\Lambda$ is icc and $[\Lambda:\Sigma_1\Sigma_2]<\infty$. Therefore, $[\Delta_2:\Sigma_2]<\infty$ and $\Delta_1$ and $\Delta_2$ are commuting subgroups of $\Lambda$.

Remark that since $[\Sigma_1:\Delta_1]<\infty$, we obtain that $\mathcal O_{\Sigma_1}(g)$ is finite if and only if $\mathcal O_{\Delta_1}(g)$ is finite, for any $g\in\Lambda$. This implies that $\{g\in\Lambda|\mathcal O_{\Delta_1}(g) \text{ is finite}\} $ equals $\Delta_2.$ Hence,
\begin{equation}\label{same}
L(\Delta_1)'\cap M=(A^{\Delta_2}\rtimes\Delta_1)'\cap M =A^{\Delta_1}\rtimes\Delta_2. 
\end{equation}

Indeed, first note that $L(\Delta_1)'\cap M \supset(A^{\Delta_2}\rtimes\Delta_1)'\cap M \supset A^{\Delta_1}\rtimes\Delta_2. $ Now, take $x\in L(\Delta_1)'\cap M.$ 
Since $\{g\in\Lambda|\mathcal O_{\Delta_1}(g) \text{ is finite}\}$ equals $\Delta_2$, it follows that $x\in A\rtimes\Delta_2.$ Using the fact that $\Delta_1$ and $\Delta_2$ commute, we obtain that $x\in A^{\Delta_1}\rtimes\Delta_2$, which proves \eqref{same}.

We continue by showing the following claim.

{\bf Claim.} We have $(A^{\Delta_2}\rtimes\Delta_1)e\prec_M P_1$ and $e(A^{\Delta_1}\rtimes\Delta_2)e\prec^s_M P_2.$

{\it Proof of the Claim.} The assumption implies that $L(\Delta_1)e\prec_M P_1.$ By passing to relative commutants and by applying twice Lemma \ref{va}(3), we obtain that $((L(\Delta_1)'\cap M)'\cap M)e\prec_M P_1.$ Relation \eqref{same} readily implies that $(A^{\Delta_2}\rtimes\Delta_1)e\prec_M P_1$.

For proving the second statement, we first show $(A^{\Delta_1}\rtimes\Sigma_2)e\prec^s_M P_2.$
Denote $\Omega_1=\{g\in\Lambda|\mathcal O_{\Sigma_2}(g)$  is finite $\}$
and let $\Omega_2=C_{\Sigma_2}(\Omega_1)$. One can show as before that the inclusions of subgroups $\Delta_1\subset\Sigma_1\subset \Omega_1$ and $\Omega_2\subset\Sigma_2\subset \Delta_2$ are of finite index. As in \eqref{same}, we obtain that $L(\Omega_2)'\cap M= A^{\Omega_2}\rtimes\Omega_1.$ Let $f\in
\mathcal N_{eMe}((A^{\Delta_1}\rtimes\Sigma_2)e)'\cap eMe\subset (L(\Delta_1\Sigma_2)e)'\cap eMe= A^{\Delta_1\Sigma_2}$e. The last equality follows from the fact that $[\Lambda:\Delta_1\Sigma_1]<\infty$ and $\Lambda$ is icc. The assumption implies $L(\Omega_2)f\prec_M P_2$. As in the previous paragraph, by applying Lemma \ref{va}(3), we get that $(L((\Omega_2)'\cap M)'\cap M)f\prec_M P_2$. This shows that $(A^{\Omega_1}\rtimes\Omega_2)f\prec_M P_2$, since $L(\Omega_2)'\cap M=A^{\Omega_2}\rtimes\Omega_1$. 

Since $[\Omega_1:\Delta_1]<\infty$ and $[\Sigma_2:\Omega_2]<\infty$, we can use Lemma \ref{fi3} and obtain that $A^{\Omega_1}\rtimes\Omega_2\subset A^{\Delta_1}\rtimes\Sigma_2$ has finite index. Note that $\Omega_1$ and $\Omega_2$ are icc since $\Lambda$ is icc and $[\Lambda:\Omega_1\Omega_2]<\infty$. This shows that $\mathcal Z(A^{\Omega_1}\rtimes\Omega_2)=A^{\Omega_1\Omega_2}$ and is completely atomic by Corollary \ref{atomic}. Hence, Lemma \ref{fi2} combined with Lemma \ref{extension} imply that $(A^{\Delta_1}\rtimes\Sigma_2)f\prec_M^{s'}(A^{\Omega_1}\rtimes\Omega_2)f.$ In combination with the conclusion of the previous paragraph and by applying Lemma \ref{transitivity}, we get that $(A^{\Delta_1}\rtimes\Sigma_2)f\prec_M P_2.$ This shows that $(A^{\Delta_1}\rtimes\Sigma_2)e\prec^s_M P_2.$

Finally, for finishing the proof of the claim, note that Lemma \ref{fi3} shows that $[A^{\Delta_1}\rtimes\Delta_2:A^{\Delta_1}\rtimes\Sigma_2]<\infty$. This gives by applying Lemma \ref{fi1} that $e(A^{\Delta_1}\rtimes\Delta_2)e\prec_M^s (A^{\Delta_1}\rtimes\Sigma_2)e$ since $e\in A^{\Delta_1}.$ Therefore, by applying Lemma \ref{transitivity}(1) and using the conclusion of the previous paragraph, we obtain that $e(A^{\Delta_1}\rtimes\Delta_2)e\prec^s_M P_2.$
\hfill$\square$

Lemma \ref{atomic} shows that $A^{\Delta_1\Delta_2}$ is completely atomic. Thus, there exists a projection $f\in A^{\Delta_1\Delta_2}$ such that $A^{\Delta_1\Delta_2}f=\mathbb C f$, $ef\neq 0$ and $(A^{\Delta_2}\rtimes\Delta_1)ef\prec_M P_1.$

Denote $t=\tau(ef)$ and $M_0=efMef$. Note that we have the identification $M_0=M^t=P_1^t\bar\otimes P_2.$ For ease of notation, we denote $Q_1=P_1^t$ and $Q_2=P_2.$

Note that $\Delta_1$ and $\Delta_2$ are icc as well and that $(A^{\Delta_2}\rtimes\Delta_1)ef$ and $e(A^{\Delta_1}\rtimes\Delta_2)ef$ are II$_1$ factors.
Thus, by combining the Claim with \eqref{same}, we have that the von Neumann algebras $(A^{\Delta_2}\rtimes\Delta_1)ef\prec_M Q_1$, and $(A^{\Delta_2}\rtimes\Delta_1)ef$ and $[(A^{\Delta_2}\rtimes\Delta_1)ef]'\cap M_0$ are II$_1$ factors.
Therefore, we can apply \cite[Proposition 12]{OP03} and deduce that there exist a decomposition $M_0=Q_1^{s}\bar\otimes Q_2^{1/s}$, for some $s>0$ and a unitary $v\in M_0$ such that 
\begin{equation}\label{eq1}
v(A^{\Delta_2}\rtimes\Delta_1)efv^*\subset Q_1^{s}.
\end{equation}

Using once again relation \eqref{same} we obtain that 
$$Q_2^{1/s}\subset ve(A^{\Delta_1}\rtimes \Delta_2)ef v^*.$$ 

By applying \cite[Theorem A]{Ge95}, we can find a factor $C\subset Q_1^{s}$ such that $C\bar\otimes Q_2^{1/s}=ve(A^{\Delta_1}\rtimes \Delta_2)ef v^*.$ Note that the Claim shows that $e(A^{\Delta_1}\rtimes \Delta_2)ef\prec_{M_0} Q_2$, which implies that $C\prec_{M_0} Q_2$. Finally, using that $C\subset Q_1^{s}$, we obtain that $C$ is not diffuse by Lemma \ref{not diffuse}. Since $C$ is a factor, it must be finite dimensional, hence $C=\mathbb M_{k}(\mathbb C),$ for some $k\ge 1.$ Denoting $t_0=s/k$, we get that 
\begin{equation}\label{z}
Q_2^{1/t_0}=ve(A^{\Delta_1}\rtimes \Delta_2) efv^*.
\end{equation}

Denote $\Omega_1=\{g\in\Lambda| \mathcal O_{\Sigma_2}(g) \text{ is finite}\}$ and $\Omega_2=C_{\Sigma_2}(\Omega_1)$, as in the proof of the Claim. Remark that $\Omega_1=\{g\in\Lambda| \mathcal O_{\Delta_2}(g) \text{ is finite}\}$ is normalized by $\Delta_2$ and $\Delta_2=\{g\in\Lambda| \mathcal O_{\Omega_1}(g) \text{ is finite}\}$ is normalized by $\Omega_1$, since $[\Delta_2:\Sigma_2]<\infty$ and $[\Omega_1:\Sigma_1]<\infty$. Note also that $\Omega_1\cap \Delta_2=1$ since $\Delta_2$ is icc. Hence, any commutator $[g,h]$, with $g\in\Omega_1,h\in\Delta_2$ belongs to $\Omega_1\cap \Delta_2=1.$ This shows that $\Omega_1$ and $\Delta_2$ are commuting subgroups.

Using \eqref{z}, we get that $Q_1^{t_0}=v((A^{\Delta_1}\rtimes\Delta_2)'\cap M)efv^*\subset vef(A^{\Omega_2}\rtimes\Omega_1)efv^*$. 
It follows that $M_0=Q_1^{t_0}\bar\otimes Q_2^{1/t_0}\subset vef(A\rtimes\Omega_1\Delta_2)efv^*$. We claim that $\Lambda_1:=\Omega_1$ and $\Lambda_2:=\Delta_2$ satisfy the conclusions of the theorem.

Note that the previous inclusion implies that $efMef=ef(A\rtimes\Lambda_1\Lambda_2)ef$. Let $z$ be the central support of $ef$ in $A\rtimes\Lambda_1\Lambda_2$ and note that $zMz=(A\rtimes\Lambda_1\Lambda_2)z.$ Denote by $Y_0\subset Y$ the $\Lambda_1\Lambda_2$-invariant measurable subset such that $z=1_{Y_0}\in\mathcal Z(A\rtimes\Lambda_1\Lambda_2)=A^{\Lambda_1\Lambda_2}$ and note that $\Lambda\car Y$ is induced from $\Lambda_1\Lambda_2\car Y_0$, since $zu_gz=0$ for any $g\in\Lambda\setminus\Lambda_1\Lambda_2$. Here, we denote by $\{u_g\}_{g\in\Lambda}$ the canonical unitaries which implement the action $\Lambda\car A$.

To this end, we use the following identification $zMz=M^{\tau(z)}=P_1^{\tau(z)}\bar\otimes P_2$. Observe that \eqref{z} shows that $L(\Lambda_2)z\prec_{M^{\tau(z)}} P_2.$ 
By passing to relative commutants and by applying twice Lemma \ref{va}(3) we obtain that $(A^{\Lambda_1}\rtimes\Lambda_2)z\prec_{M^{\tau(z)}} P_2.$ Therefore, by applying \cite[Proposition 12]{OP03} (note that $(A^{\Lambda_1}\rtimes\Lambda_2)z$ and $((A^{\Lambda_1}\rtimes\Lambda_2)z)'\cap zMz$ are II$_1$ factors) and \cite[Theorem A]{Ge95} as before, there exist a decomposition $M^{\tau(z)}=P_1^{t_1}\bar\otimes P_2^{t_2}$, for some $t_1,t_2>0$ with $t_1t_2=\tau(z)$, a unitary $w\in M^{\tau(z)}$, and a von Neumann subalgebra $D\subset P_2^{t_2}$ such that $P_1^{t_1}\bar\otimes D=w(A^{\Lambda_2}\rtimes\Lambda_1)zw^*$. Remark that Lemma \ref{fi3} and Lemma \ref{fi1} show that $[(A^{\Lambda_2}\rtimes\Lambda_1)z:(A^{\Delta_2}\rtimes\Delta_1)z]<\infty$, and Lemma \ref{fi2} implies $(A^{\Lambda_2}\rtimes\Lambda_1)z\prec_{(A^{\Lambda_2}\rtimes\Lambda_1)z}^{s'}(A^{\Delta_2}\rtimes\Delta_1)z$. Recall that the Claim shows that $(A^{\Delta_2}\rtimes\Delta_1)z\prec_M P_1$, since $ef\leq z$. Now, Lemma \ref{extension} and Lemma \ref{transitivity}(2) imply that $(A^{\Lambda_2}\rtimes\Lambda_1)z\prec_{M^{\tau(z)}} P_1^{t_1}$, hence, $D\prec_{M^{\tau(z)}} P_1^{t_1}$. Since $D\subset P_2^{t_2}$, it cannot be diffuse by Lemma \ref{not diffuse}. Using that $D$ is factor, it follows that it must be finite dimensional. Therefore, $D=\mathbb M_{p}(\mathbb C)$, for some $p\ge 1$. Thus, $P_1^{t_1p}=w(A^{\Lambda_2}\rtimes\Lambda_1)zw^*$. By passing to relative commutants, we obtain that $P_2^{t_2/p}=w(A^{\Lambda_1}\rtimes\Lambda_2)zw^*$. This clearly implies the conclusion of the theorem by noticing that $zMz=L^\infty(Y_0)\rtimes\Lambda_1\Lambda_2$ and by representing $A^{\Lambda_2}z=L^\infty(Y_1)$ and $A^{\Lambda_1}z=L^\infty(Y_2)$ for some standard probability spaces $(Y_1,\nu_1)$ and $(Y_2,\nu_2)$. 
\hfill$\blacksquare$

\section{Proof of Theorem \ref{A2}}

In this section we will prove the von Neumann algebraic version of Theorem \ref{A2} and use it to derive the main result of the introduction. First, we make the remark that if $\Lambda\car Y$ is induced from an action $\Lambda_0\car Y_0$, then we have the identification $L^\infty(Y_0)\rtimes\Lambda_0^{[\Lambda:\Lambda_0]}=L^\infty(Y)\rtimes\Lambda$ since $L^\infty(Y_0)\rtimes\Lambda_0=1_{Y_0}(L^\infty(Y)\rtimes\Lambda)1_{Y_0}$. 

\begin{theorem}\label{vNa}
Let $\Gamma_1,\dots,\Gamma_n$ be property (T), biexact, weakly amenable groups and denote $\Gamma=\Gamma_1\times ... \times\Gamma_n$. For every $1\leq i\leq n$, let $\Gamma_i\overset{}{\car} (X_i,\mu_i)$ be a free ergodic pmp action and 
denote $M_i=L^\infty (X_i)\rtimes_{}\Gamma_i $ and $M=M_1\bar\otimes \dots \bar\otimes M_n.$

Let $\Lambda\car (Y,\nu)$ be a free ergodic pmp action of a countable icc group $\Lambda$ such that $M^t=L^\infty(Y)\rtimes\Lambda,$ for some $t>0.$

Then $\Lambda\car Y$ is induced from an action $\Lambda_0\car Y_0$ and there exist a decomposition $\Lambda_0=\Lambda_1\times ... \times \Lambda_n$ and pmp actions $\Lambda_i\car Y_i$, positive numbers $t_i>0$ with $t_1\dots t_n=t/[\Lambda:\Lambda_0]$ and a unitary $u\in M^{t/[\Lambda:\Lambda_0]}$ such that $\Lambda_0\car Y_0$ is isomorphic to the product action $\Lambda_1\times ... \times \Lambda_n\car Y_1\times ...\times Y_n$ and $u(L^\infty(Y_i)\rtimes\Lambda_i) u^*=M_i^{t_i},$ for every $1\leq i \leq n$.

\end{theorem}

An essential ingredient of the proof of Theorem \ref{vNa} consists of applying A. Ioana's ultrapower technique \cite{Io11}, which we recall in the following form. This result is essentially contained in the proof of \cite[Theorem 3.1]{Io11} and its statement is roughly \cite[Theorem 4.1]{DHI16}. We leave the proof to the reader, since it follows verbatim the proof of the result in \cite{DHI16}.

\begin{theorem}[\!\!{\cite{Io11}}]\label{Th: ultrapower}
Let $M=B\rtimes\Lambda$ be a II$_1$ factor, where $\Lambda\car B$ is a trace preserving action on a tracial von Neumann algebra. Let $\Delta:M\to M\bar\otimes M$ be the $*$-homomorphism defined by $\Delta(b)=b\otimes 1$ and $\Delta (bv_\lambda)=bv_\lambda\otimes v_\lambda$, for all $b\in B$ and $\lambda\in\Lambda.$ Let $P,Q\subset M$ be von Neumann subalgebras such that $\Delta(P)\prec_{M\bar\otimes M}M\bar\otimes Q$.

Then there exists a decreasing sequence of subgroups $\Sigma_k<\Lambda$ such that $P\prec_M B\rtimes\Sigma_k$, for every $k\ge 1$, and $Q'\cap M\prec_M B\rtimes (\cup_{k\ge 1} C_\Lambda(\Sigma_k)).$
\end{theorem}

Note that the ultrapower technique has recently been used in several other works \cite{CdSS15, KV15, DHI16, CI17, CU18}.

{\it Proof of Theorem \ref{vNa}.}
We start the proof by fixing some notation. Let $l\ge t$ be an integer and let $p\in L^\infty(X\times \mathbb Z/l\mathbb Z)$ be a projection of trace $t/l$ such that $M^t=p(L^\infty(X\times \mathbb Z/l\mathbb Z)\rtimes(\Gamma\times\mathbb Z/l\mathbb Z))p$. For ease of notation, we assume that $0<t\leq 1$, and therefore we can take $l=1$. Hence, $pMp=L^\infty(Y)\rtimes\Lambda.$

For any $1\leq i \leq n$, let $\hat i=\{1,\dots,n\}\setminus \{i\}$ and $A_i=L^\infty(X_i)$. 
For a subset $F\subset \{1,\dots,n\}$, we denote $\Gamma_F=\times _{i\in F}\Gamma_i$, $A_F=\bar\otimes_{i\in F} A_i$ and $M_F=\bar\otimes_{i\in F}M_i.$ 
Since the groups $\Gamma_i$'s are weakly amenable and biexact, \cite[Theorem 1.3]{PV12} implies that $M$ has a unique Cartan subalgebra up to unitarily conjugacy. Hence, we may assume that $L^\infty(Y)=L^\infty(X)p.$ Denote $A=L^\infty(X)$ and $B=L^\infty(Y).$

Following \cite{PV09} we define the comultiplication $\Delta: M^t\to M^t\bar\otimes L(\Lambda)$ by $\Delta(bv_\lambda)=bv_\lambda\otimes v_\lambda$, for all $b\in L^\infty(Y)$ and $\lambda\in \Lambda.$ 
Since $M=M_{\hat j}\bar\otimes M_j$, we can write $M^t=M_{\hat j}\bar\otimes M_j^t$, for every $1\leq j\leq n.$
The proof of the theorem is divided between the following four claims.

{\bf Claim 1.} We can find $1\leq j_0\leq n$ such that $\Delta(L(\Gamma_{\hat n}))\prec_{M^t\bar\otimes M^t} M^t\bar\otimes M_{\hat j_0}.$

{\it Proof of Claim 1.} First, remark that there exists $1\leq j_0 \leq n$ such that $\Delta(M_n^t)$ is not amenable relative to $M^t\bar\otimes M_{\hat j_0}$ inside $M^t\bar\otimes M^t$. Otherwise, by applying Proposition \ref{intersection}(2), we obtain that $\Delta(M^t_n)$ is amenable relative to $M^t\otimes 1$ inside $M^t\bar\otimes M^t$. This implies by \cite[Lemma 10.2]{IPV10} that $M_n$ is amenable, contradiction. Therefore, by using the fact that $\Gamma_{j_0}\in \mathcal C_{rss}$, we obtain that $\Delta(M_{\hat n})\prec_{M^t\bar\otimes M} M^t\bar\otimes (A_{j_0}\bar\otimes M_{\hat j_0}).$ Since $\Gamma_{\hat n}$ has property (T), Lemma \ref{L: rigid} shows that $\Delta(L(\Gamma_{\hat n}))\prec_{M^t\bar\otimes M^t} M^t\bar\otimes M_{\hat j_0}.$  
\hfill$\square$

We are now in a position to apply the ultrapower technique from \cite{Io11}. Combining Claim 1 with Theorem \ref{Th: ultrapower}, we deduce the existence of a decreasing sequence of subgroups $\Sigma_k<\Lambda$ such that 
\begin{equation}\label{ultrapower}
L(\Gamma_{\hat n})\prec_{M^t} B\rtimes\Sigma_k,  \text{ for every } k\ge 1, \text{ and }M_{j_0}^t\prec_{M^t} B\rtimes (\cup_{k\ge 1}C_\Lambda(\Sigma_k)).
\end{equation}

{\bf Claim 2.} There exists a non-amenable subgroup $\Sigma< \Lambda$ with non-amenable centralizer $C_\Lambda(\Sigma)$ and a projection $e\in \mathcal N_{pMp}(B\rtimes\Sigma)'\cap pMp$ such that 
$$
A\rtimes\Gamma_{\hat n} \prec_M^{s'} (B\rtimes\Sigma) e \text{   and    }( B\rtimes\Sigma) e\prec_M^s A\rtimes\Gamma_{\hat n}.
$$

{\it Proof of Claim 2.} Relation \eqref{ultrapower} implies that there exists a non-amenable subgroup $\Sigma< \Lambda$ with non-amenable centralizer $C_\Lambda(\Sigma)$ such that $L(\Gamma_{\hat n})\prec_M B\rtimes\Sigma$. We can use Lemma \ref{union}(1) and derive that $A\rtimes\Gamma_{\hat n} \prec_M B\rtimes\Sigma$. By applying Lemma \ref{dhi}(2), there exists a non-zero projection $e\in \mathcal N_{pMp}(B\rtimes\Sigma)'\cap pMp$ such that $A\rtimes\Gamma_{\hat n} \prec_M^{s'} (B\rtimes\Sigma) e$. 

For proving the claim, it remains to show that $( B\rtimes\Sigma) e\prec_M^s A\rtimes\Gamma_{\hat n}$. 
To this end, take a projection $f\in \mathcal N_{pMp}(B\rtimes\Sigma)'\cap pMp\subset B^{\Sigma C_\Lambda(\Sigma)}$ with $f\leq e.$ First, note that $L(C_{\Lambda}(\Sigma))f$ is amenable relative to $A\rtimes\Gamma_n.$ Indeed, take an arbitrary $1\leq i \leq n-1$. Since $\Gamma_i$ belongs to $\mathcal C_{rss}$, Lemma \ref{L: rss} implies that $L(\Sigma)f\prec_M A\rtimes\Gamma_{\hat i}$ or $L(C_\Lambda(\Sigma))f$ is amenable relative to $A\rtimes\Gamma_{\hat i}$. If the former relation holds, we can apply Lemma \ref{union}(1) and obtain that $(B\rtimes\Sigma )f\prec_M A\rtimes\Gamma_{\hat i}$.
In combination with $A\rtimes\Gamma_{\hat n}\prec_M^{s'} (B\rtimes\Sigma)f$, Lemma \ref{transitivity}(2) implies that $A\rtimes\Gamma_{\hat n}\prec_M A\rtimes\Gamma_{\hat i}$, contradicting the fact that $\Gamma_i$ is an infinite group. 
Therefore, $L(C_\Lambda(\Sigma))f$ is amenable relative to $A\rtimes\Gamma_{\hat i}$, for all $1\leq i\leq n-1$. By applying Proposition \ref{intersection}(2), we obtain that $L(C_{\Lambda}(\Sigma))f$ is amenable relative to $A\rtimes\Gamma_n.$

Since $\Gamma_n$ belongs to $\mathcal C_{rss}$, we apply once again Lemma \ref{L: rss} and obtain that $L(C_{\Lambda}(\Sigma))f$ is amenable relative to $A\rtimes\Gamma_{\hat  n}$ or $L(\Sigma)f\prec_M A\rtimes\Gamma_{\hat n}$. The former relation combined with the fact that $L(C_{\Lambda}(\Sigma))f$ is amenable relative to $A\rtimes\Gamma_n$ gives that $C_\Lambda(\Sigma)$ is amenable by Proposition \ref{intersection}(2), contradiction. Hence, $L(\Sigma)f\prec_M A\rtimes\Gamma_{\hat n}$. By applying Lemma \ref{union}(1), we obtain that $(B\rtimes\Sigma) f\prec_M A\rtimes\Gamma_{\hat n}.$ This shows that $(B\rtimes\Sigma)e\prec_M^s A\rtimes\Gamma_{\hat n},$ which ends the proof of the claim.
\hfill$\square$

Remark that Claim 2 allows us to apply Proposition \ref{ME} and derive that $\Sigma$ is measure equivalent to $\Gamma_{\hat n}.$ Since property (T) is a measure equivalence invariant \cite[Corollary 1.4]{Fu99a}, we deduce that $\Sigma$ has property (T) as well.

Denote $\Delta=\{g\in\Lambda|\mathcal O_\Sigma(g) \text{ is finite}\}$. Note that $\Delta$ is normalized by $\Sigma$ and $L(\Sigma)'\cap (B\rtimes\Lambda)\subset B\rtimes \Delta.$

{\bf Claim 3.} $B^{\Delta\Sigma}e$ is completely atomic and $\Delta\Sigma$ is a finite index subgroup of $\Lambda$. 

{\it Proof of Claim 3.} First, we show that $B^{\Delta\Sigma}e$ is completely atomic and use this to derive the second part of the claim. Claim 2 implies that $A\rtimes\Gamma_{\hat n}\prec_M (B\rtimes \Delta\Sigma)f,$ for all $f\in B^{\Delta\Sigma}$ such that $f\leq e.$ By passing to relative commutants and by applying Lemma \ref{va}(3), we obtain that $B^{\Delta\Sigma}f\prec_M A_n.$ Notice that $\mathcal N_{eMe}(B^{\Delta\Sigma}e)'\cap eMe\subset B^{\Delta\Sigma} e$, and hence, Lemma \ref{dhi}(1) gives that $B^{\Delta\Sigma}e\prec_M^s A_n.$ 

Now we show that $B^{\Delta\Sigma}e\prec_M^s M_{\hat n}.$ Take a non-zero projection $f\in B^{\Delta\Sigma}$ with $f\leq e.$ Claim 2 implies that $L(\Sigma)f\prec_M A\rtimes\Gamma_{\hat n}$. Since $\Sigma$ has property (T), by using Lemma \ref{L: rigid} we deduce that $L(\Sigma)f\prec_M M_{\hat n}.$ By passing to relative commutants and by applying Lemma \ref{va}(3), it follows that $((L(\Sigma)'\cap pMp)'\cap pMp) f\prec_M M_{\hat n}.$ Since $L(\Sigma)'\cap pMp\subset B\rtimes\Delta,$ we get $B^{\Sigma\Delta}f\prec_M M_{\hat n}.$ Therefore, 
$B^{\Delta\Sigma}e\prec_M^s M_{\hat n}.$
Together with the conclusion of the previous paragraph,  Proposition \ref{intersection}(1) shows that $B^{\Delta\Sigma}e\prec_M^s \mathbb C 1$, which implies that $B^{\Delta\Sigma}e$ is completely atomic.  

We continue by taking a non-zero projection $e_0\in B^{\Delta\Sigma}$ with $e_0\leq e$ such that $B^{\Delta\Sigma}e_0=\mathbb C e_0.$ On one hand, by combining Lemma \ref{L: rigid} and Claim 2 it follows that $L(\Sigma)e_0\prec_M M_{\hat n}.$ By considering relative commutants Lemma \ref{va}(3) implies that $M_n\prec_M (B\rtimes\Delta\Sigma)e_0.$
On the other hand, Claim 2 gives that $M_{\hat n}\prec_M (B\rtimes\Delta\Sigma)e_0$. 
Note that $\mathcal N_{e_0Me_0}((B\rtimes\Delta\Sigma)e_0)'\cap e_0Me_0\subset B^{\Delta\Sigma}e_0=\mathbb C e_0.$ Therefore, by applying Lemma \ref{L: full}, we obtain that $M\prec_M (B\rtimes\Delta\Sigma)e_0$. It is easy to see that this implies $L(\Lambda)\prec_{L(\Lambda)} L(\Delta\Sigma).$ Hence, we obtain that $\Delta\Sigma$ is a finite index subgroup of $\Lambda$ by \cite[Lemma 2.5(1)]{DHI16}.
\hfill$\square$

{\bf Claim 4.} There exists a finite index subgroup  $\Sigma_0$ of $\Sigma$ such that $[\Lambda:\Sigma_0 C_\Lambda(\Sigma_0)]<\infty$. Moreover, there exists a non-zero projection $q\in B^{\Sigma_0C_\Lambda(\Sigma_0)}$ satisfying
$$
L(\Sigma_0)q\prec_{M^t} M_{\hat n} \text{ and  } L(C_\Lambda(\Sigma_0))q\prec^s_{M^t} M^t_n.
$$
{\it Proof of Claim 4.} Recall that $\Delta=\{g\in\Lambda| \mathcal O_\Sigma(g) \text{ is finite}\}.$ Let $\{\mathcal O_k\}_{k\in\mathbb N}$ be a countable enumeration of all the finite orbits of the action by conjugation of $\Sigma$ on $\Lambda$ and notice that $\Delta=\cup_{k\in\mathbb N}\mathcal O_k.$ Denote $\mathcal S_k=\cup_{i=1}^k \mathcal O_i$ and note that $\Delta_k:=\langle \mathcal S_k \rangle$ is an ascending sequence of subgroups of $\Lambda$ normalized by $\Sigma$ satisfying $\cup_k \Delta_k=\Delta$. 

Since $\Lambda$ is measure equivalent to a property (T) group and $[\Lambda: \Delta\Sigma]<\infty$, we derive by \cite[Corollary 1.4]{Fu99a} that $\Delta\Sigma$ has property (T) as well. Therefore, there exists $k\in\mathbb N$ such that $\Delta\Sigma=\Delta_k\Sigma.$ Since $\mathcal S_k$ is a finite set, then the subgroup $\Sigma_0:=\cap_{g\in\mathcal S_k}C_\Sigma (g)$ has finite index in $\Sigma$ and commutes with $\Delta_k$. Thus, $[\Lambda:\Delta_k \Sigma_0]<\infty$, which shows that $[\Lambda:\Sigma_0 C_\Lambda(\Sigma_0)]<\infty$.

For proving the moreover part, remark first that as in Claim 2
we can find a non-zero projection
$q\in \mathcal N_{pMp}(B\rtimes\Sigma_0)'\cap pMp\subset B^{\Sigma_0C_\Lambda(\Sigma_0)}$
such that
\begin{equation}\label{done}
A\rtimes\Gamma_{\hat n}\prec^{s'}_M (B\rtimes\Sigma_{0})q  \text{ and }  (B\rtimes\Sigma_0)q \prec^s_M A\rtimes\Gamma_{\hat n}.
\end{equation}

We continue by showing that $(B\rtimes C_\Lambda(\Sigma_0))q\prec^s_{M} A\rtimes\Gamma_n.$ Let $f \in \mathcal N_{pMp}(B\rtimes C_\Lambda(\Sigma_0))'\cap pMp\subset B^{\Sigma_0C_\Lambda(\Sigma_0)}$ such that $f\leq q.$
First, note that $L(\Sigma_0)f\nprec_M A\rtimes\Gamma_{\hat i}$, for any $1\leq i\leq n-1.$ Indeed, if there exists such an $i$ for which $L(\Sigma_0)f\prec_M A\rtimes\Gamma_{\hat i}$, we get by Lemma \ref{union}(1) that $(B\rtimes\Sigma_0)f\prec_M A\rtimes\Gamma_{\hat i}$. By using Lemma \ref{transitivity}(2), relation \eqref{done} shows that $\Gamma_i$ is a finite group, contradiction.

Note that actually $L(\Sigma_{0})f$ is not amenable relative to $A\rtimes\Gamma_{\hat i}$, for any $1\leq i\leq n-1.$ Indeed, if there exists such an $i$, by using that $\Gamma_i\in \mathcal C_{rss}, $ we get that $L(\Sigma_0C_\Lambda(\Sigma_0))f$ is amenable relative to $A\rtimes\Gamma_{\hat i}$ since $L(\Sigma_0)f\nprec_M A\rtimes\Gamma_{\hat i}$. By using Lemma \ref{union}(2) we obtain that $B\rtimes(\Sigma_0C_\Lambda(\Sigma_0))f_0$ is amenable relative to $A\rtimes\Gamma_{\hat i},$ for some non-zero projection $f_0\in B^{\Sigma_0C_\Lambda(\Sigma_0)}$. By combining Lemma \ref{fi3}, Lemma \ref{fi1} and \cite[Lemma 2.6(3)]{DHI16}, it follows that $f_0(B\rtimes\Lambda)f_0$ is amenable relative to $B\rtimes(\Sigma_0C_\Lambda(\Sigma_0))f_0$ inside $M$. 
By applying \cite[Proposition 2.4(3)]{OP07} and \cite[Lemma 2.6(2)]{DHI16} we get that $A\rtimes\Gamma$ is amenable relative to $A\rtimes\Gamma_{\hat i}$. Hence, \cite[Proposition 2.4(1)]{OP07} shows that $\Gamma_i$ is amenable, false.

Now, we are finally showing that $(B\rtimes C_{\Lambda}(\Sigma_0))f\prec^s_M A\rtimes\Gamma_{n}.$ Let $1\leq i\leq n-1.$ By using once again that $\Gamma_i\in \mathcal C_{rss}$ we deduce by Lemma \ref{L: rss} that $L(C_\Lambda(\Sigma_0))f\prec_M A\rtimes\Gamma_{\hat i}$ since $L(\Sigma_0)f$ is not amenable relative to $A\rtimes\Gamma_{\hat i}$ inside $M$. Note that $(B\rtimes C_\Lambda(\Sigma_0))f\prec_M A\rtimes\Gamma_{\hat i}$ by Lemma \ref{union}(1), and hence $(B\rtimes C_\Lambda(\Sigma_0))q\prec_M^{s} A\rtimes\Gamma_{\hat i}$ by Lemma \ref{dhi}(1).
We can apply now Proposition \ref{intersection}(1) and deduce that $(B\rtimes C_{\Lambda}(\Sigma_0))q\prec^s_ M A\rtimes\Gamma_{n}.$

By applying Lemma \ref{transitivity}(1), we derive that $L(\Sigma_0)q\prec^s_{M} A\rtimes\Gamma_{\hat n} $ and $L(C_{\Lambda}(\Sigma_0))q\prec^s_M A\rtimes\Gamma_n$.
Since the groups $\Sigma_0$ and $C_{\Lambda}(\Sigma_0)$ have property (T), we can use Lemma \ref{L: rigid} and obtain the claim.
\hfill$\square$

Finally, by applying Theorem \ref{P: splits} we obtain that there exist commuting subgroups $\Lambda_{1}^{n-1},\Lambda_n<\Lambda$ such that $\Lambda\car Y$ is induced from an action $\Lambda_{1}^{n-1}\Lambda_{n}\car Y_n^n$. Denote $p_n=1_{Y^n_n}\in B^{\Lambda_{1}^{n-1}\Lambda_{n}}.$
Moreover, there exist a decomposition $M^{\tau(p_n)}=M_{\hat n}^{t_1^{n-1}}\bar\otimes M_{n}^{t_n}$ for some $t_1^{n-1},t_n>0$ with $t_1^{n-1}t_n=\tau(p_n)$, and a unitary $u_n\in M^{\tau(p_n)}$ such that:
$$
u_n(B^{\Lambda_n}\rtimes\Lambda_1^{n-1})p_nu_n^*=M_{\hat n}^{t_1^{n-1}} \text{ and } u_n(B^{\Lambda_1^{n-1}}\rtimes\Lambda_n)p_nu_{n}^*=M_n^{t_n}.
$$
In particular, there exist pmp actions $\Lambda_{1}^{n-1}\Lambda_{n}\car Y_1^{n-1}$ and $\Lambda_n\car Y_n$ such that 
\begin{equation}\label{induction1}
u_n(L^\infty(Y_1^{n-1})\rtimes\Lambda_1^{n-1})u_n^*=M_{\hat n}^{t_1^{n-1}} \text{ and } u_n(L^\infty(Y_n)\rtimes\Lambda_n)u_n^*=M_n^{t_n},
\end{equation}
and
$\Lambda_{1}^{n-1}\Lambda_{n}\car Y_n^n$ is isomorphic to the product action $\Lambda_{1}^{n-1}\times \Lambda_n\car Y_1^{n-1}\times Y_n$.

Applying an induction argument and Remark \ref{induce} to relation \eqref{induction1}, it is easy to see that the conclusion follows. 
\hfill$\blacksquare$





Before we proceed with the proof of Theorem \ref{A2}, we recall some notation. Assume that $M$ is a II$_1$ factor and $A\subset M$ is a {\it Cartan subalgebra}, i.e. a maximal abelian regular von Neumann subalgebra. The inclusion $A^t\subset M^t$ is defined as the isomorphism class of the inclusion $p(\ell^\infty(\mathbb Z)\bar\otimes A)p\subset p(\mathbb B(\ell^2(\mathbb Z))\bar\otimes M)p$, where $p\in \mathbb B(\ell^2(\mathbb Z))\bar\otimes M$ is a projection satisfying  $(\text{Tr}\otimes\tau)(p)=t$. Here, we denote by $\ell^\infty(\mathbb Z)\subset\mathbb B(\ell^2(\mathbb Z))$ the subalgebra of diagonal operators and by $\text{Tr}$ the usual trace on $\mathbb B(\ell^2(\mathbb Z))$.

{\it Proof of Theorem \ref{A2}.}
For any $1\leq i\leq n$, denote $M_i=L^\infty(X_i)\rtimes\Gamma_i.$ Applying Theorem \ref{vNa}, we obtain that that $\Lambda\car Y$ is induced from an action $\Lambda_0\car Y_0$ and there exist a decomposition $\Lambda_0=\Lambda_1\times ... \times \Lambda_n$ and pmp actions $\Lambda_i\car Y_i$, positive numbers $t_i>0$ with $t_1\dots t_n=t/[\Lambda:\Lambda_0]$ and a unitary $u\in M^{t/[\Lambda:\Lambda_0]}$ such that $\Lambda_0\car Y_0$ is isomorphic to $\Lambda_1\times ... \times \Lambda_n\car Y_1\times ...\times Y_n$ and $u(L^\infty(Y_i)\rtimes\Lambda_i) u^*=M_i^{t_i},$ for all $1\leq i \leq n$.
Note that $L^\infty(Y_i)\prec_M L^\infty (X)$ and $uL^\infty(Y_i)u^*\subset M_i^{t_i}$, for any $1\leq i\leq n$.  
This implies that $uL^\infty(Y_i)u^*\prec_{M_i^{t_i}} L^\infty(X_i)^{t_i}$. Remark that $uL^\infty(Y_i)u^*$ and $ L^\infty(X_i)^{t_i}$ are Cartan subalgebras of $M_i^{t_i}$.
Hence, by applying \cite[Theorem A.1]{Po01} and \cite[Proposition in p.829]{Po01}, we get that $uL^\infty(Y_i)u^*$ is unitarily conjugate to $L^\infty(X_i)^{t_i}$ inside $M_i^{t_i}$. Using \cite{FM75}, it follows that $\Lambda_i\car Y_i$ is stably orbit equivalent to $\Gamma_i\car X_i$ with index $t_i$.
\hfill$\blacksquare$

\end{document}